\documentclass[12pt,letterpaper,reqno]{amsart}

\usepackage{times}
\usepackage[latin1]{inputenc}
\usepackage[T1]{fontenc}
\usepackage{mathrsfs}
\usepackage{latexsym}
\usepackage[dvips]{graphics}
\usepackage{epsfig}
\usepackage{enumerate}
\usepackage{amsmath,amsfonts,amsthm,amssymb,amscd}
\input amssym.def
\input amssym.tex

\def\M{\mathbb{N}_0}
\def\N{\mathbb{N}}
\def\RM{\mathcal{M}}
\def\B{\mathcal{B}}

\newcommand\T{\rule{0pt}{2.6ex}}
\newcommand\be{\begin{equation}}
\newcommand\ee{\end{equation}}
\theoremstyle{definition}
\newtheorem{Thm}{Theorem}[section]
\newtheorem{Prop}[Thm]{Proposition}
\newtheorem{Lemma}[Thm]{Lemma}
\newtheorem{Cor}[Thm]{Corollary}

\newtheorem{Conj}[Thm]{Conjecture}
\newtheorem{Ex}{Example}
\newtheorem{Rem}{Remark}
\email{urban.larsson@chalmers.se, hegarty@chalmers.se}
\address{Mathematical Sciences, Chalmers University of Technology
and University of Gothenburg, G\"oteborg, Sweden}
\email{ aviezri.fraenkel@weizmann.ac.il}
\address{Department of Computer Science and Applied Math,
Weizmann Institute of Science, 76100 Rehovot, Israel}

\begin{document}
\title{Invariant and Dual subtraction games resolving the Duch\^ene-Rigo
Conjecture.}
\vskip 5pt
\author{ Urban Larsson, Peter Hegarty, Aviezri S. Fraenkel }
\begin{abstract}
We prove a recent conjecture of Duch\^ene and Rigo,
stating that every complementary
pair of homogeneous Beatty sequences represents the solution to
an \emph{invariant} impartial game. Here
invariance means that each available move in  a game can be played
anywhere inside the game-board. In fact, we establish such a result for a
wider class of pairs of complementary sequences, and in the
process generalize the notion of a \emph{subtraction game}.
Given a pair of complementary sequences $(a_n)$ and $(b_n)$ of
positive integers, we define a game $G$
by setting $\{\{a_n, b_n\}\}$ as invariant moves.
We then introduce the invariant game $G^\star $, whose moves are
all non-zero $P$-positions of $G$. Provided the set of non-zero
$P$-positions of $G^\star$ equals $\{\{a_n,b_n\}\}$,
this \emph{is} the desired invariant game. We give sufficient
conditions on the initial pair of sequences for this 'duality' to hold.
\end{abstract}
\maketitle
\vskip 5pt

\setcounter{equation}{0}
\section{Notation, Terminology and Statement of Results}
This note concerns 2-person, impartial games (see \cite{BCG})
played under normal (as against mis\`{e}re) rules. Let $\N$, $\M$ denote
the positive and the non-negative integers respectively. For $k\in \N$, let
$\B = \B(k) := (\M^k, \oplus , \preceq)$ denote the partially-ordered semigroup
consisting of all ordered $k$-tuples of non-negative integers, where for
elements $\boldsymbol x = (x_1, \ldots ,x_k), \boldsymbol y =
(y_1, \ldots ,y_k)$ of $\B$ one defines
\be
\boldsymbol x \oplus \boldsymbol y := (x_1 + y_1, \ldots , x_k + y_k)
\ee
and
\be
\boldsymbol x \preceq \boldsymbol y \; \Leftrightarrow \; x_i \le y_i, \;
i = 1,\ldots , k.
\ee

Hence $\boldsymbol x \prec \boldsymbol y$ if $\boldsymbol x \preceq
\boldsymbol y$ and $x_i<y_i$ for some $i$. For
$\boldsymbol y \preceq \boldsymbol x$ we define
$$\boldsymbol x\ominus\boldsymbol y:= (x_1 - y_1, \ldots , x_k - y_k).$$

We call $\B$ the {\em game board}. Let $G = G(F, \B)$ denote a game, where
for all $\boldsymbol x\in \B$,
$F(\boldsymbol x)\subset \B$ defines the set of \emph{options}
of $\boldsymbol x$ in the sense that $\boldsymbol y \in F(\boldsymbol x)$
if and only if there is a move from $\boldsymbol x$ to $\boldsymbol y$.
Formally, the {\em move\/} from $\boldsymbol x$ to $\boldsymbol y$ is
the ordered pair $(\boldsymbol x, \boldsymbol y)$. In this paper, the phrase
$\lq \boldsymbol x \rightarrow \boldsymbol y$ is an option' will often be used
synonymously with $\lq \boldsymbol y \in F(\boldsymbol x)$', in order to avoid
cumbersome notation.

Given this setting, the two players only need to (randomly) pick a
starting position
$\boldsymbol x\in \B$ and decide who plays first. Then they play
by alternating in
choosing options from $F(\cdot)$ (and moving accordingly). Although we
have announced that the last player to move wins (normal play), without
some additional assumptions there is no guarantee that the game will
terminate.

By a \emph{$k$-pile subtraction game}\footnote{Our subtraction games are
generalizations of the Nim-type subtraction games defined in \cite{BCG}.
There are some alternative names for our games that can be found in the
literature, such as \emph{Take-away games, Removal games}. By our choice
we emphasize the natural additive structure on $\B$.} we mean a
game played on $\B$ such that, for each $\boldsymbol x \in \B$,
the set $F(\boldsymbol x)\subset \B$ has the property that
$\boldsymbol y \in F(\boldsymbol x)\Rightarrow \boldsymbol y \prec
\boldsymbol x$. In the setting of invariant games (to be defined below),
it will be convenient to abuse notation and also
refer to the $k$-tuple
$\boldsymbol x \ominus \boldsymbol y \succ \boldsymbol 0$ as a {\em move}.
Observe that both options and moves are then elements of $\B$, but with
different meanings.

In this paper, whenever we refer to a \emph{(subtraction) game}
we intend a $k$-pile subtraction game. Let $G$ be a game. Then
$\mathcal{T} =
\mathcal{T}(G) := \{\boldsymbol x \mid F(\boldsymbol x) = \emptyset \}$
denotes the set of \emph{terminal} positions. Clearly
$\boldsymbol 0:=(0,\ldots,0)\in \mathcal{T}$ and $\boldsymbol 0$
is unique. It is natural
to require that $\mathcal{T}$ be a \emph{lower ideal} in the poset, that is,
if $\boldsymbol x \in \mathcal{T}$ and $\boldsymbol y \prec \boldsymbol x$,
then $\boldsymbol y \in \mathcal{T}$. Clearly, in this setting, any
game must terminate within a finite number of moves and the \emph{winner}
 is the player who makes the last move. The opponent is the {\em loser}.

Recently, Duch\^ene and Rigo \cite{DR} introduced the notion of an invariant
game. A $k$-pile subtraction game $G$ is said to be \emph{invariant} if,
for all $\boldsymbol x, \boldsymbol y\in \B$ and
$\boldsymbol r \in \B\setminus \{\boldsymbol 0 \}$,
$(\boldsymbol x \oplus \boldsymbol r)\rightarrow \boldsymbol x$
is an option whenever $(\boldsymbol y \oplus \boldsymbol r)
\rightarrow \boldsymbol y$ is. Then, in particular,
 $\boldsymbol y=\boldsymbol 0$ is an option of
$\boldsymbol r$. Therefore for invariant games we refer to the set

\be
\RM (G) := \{\boldsymbol r \in \B \mid \boldsymbol 0  \in F(\boldsymbol r) \}
\ee
as the set of all {\em (invariant) moves}{\footnote{This notation
and terminology is consistent with that employed in \cite{DR}.}} of $G$.
Knowledge of this set gives a complete description of the rules
of an invariant game. If a game is not invariant it is \emph{variant}.

A position (a game) is $P$ if all of its options are $N$. Otherwise it is $N$.
This means that the first player to move wins if and only if the game is $N$.
As usual, we shall denote by $\mathcal{P}(G)$ (resp. $\mathcal{N}(G)$)
the collection of $P$- (resp. $N$-) positions of $G$.

Finally, if $G$ is a (not necessarily invariant) game, then we can
define an invariant game $G^\star$ on the same game board by setting
\be\label{dual}
\RM (G^\star) := \mathcal{P}(G) \backslash \{\boldsymbol 0 \}.
\ee

\begin{Ex}
Define $G$ by $\RM (G) = \emptyset $. Then $\mathcal{P}(G) = \B$ and so
$\RM (G^\star ) = \B \setminus \{\boldsymbol 0\}$. This
gives $\mathcal{P}(G^\star) = \{ \boldsymbol 0 \}$, so that
in fact $\mathcal{N}(G^\star) = \RM (G^\star )$. This latter equality
does not hold in general.
For example, let $G$ rather denote 2-pile Nim. Then\footnote{A subset $R$
of $\B = \M \times \M$ is {\em symmetric} if
$(x,y) \in R \Leftrightarrow (y,x) \in R.$
(We dispense with the obvious generalisation to $k > 2$ piles.)
If the sets $\mathcal{M}(G)$ and $\mathcal{T}(G)$
are symmetric subsets of $\B$, then so are the sets $\mathcal{N}(G)$
and $\mathcal{P}(G)$. In this case the game $G$ will be called
\emph{symmetric}.
Sometimes it will be convenient to denote
moves and positions of a symmetric game by unordered pairs $\{r, s\}$.
Hence, whenever we write $\lq \{r,s\} \in \RM (G)$' for example,
what we mean is that $\{(r,s), (s,r)\} \subseteq \RM (G)$.}
$\RM (G) = \{ \{0,x\}\mid x\in \N \}$ and
$\mathcal{P}(G) = \{ \{x,x \}\mid x\in \M \}$. By (\ref{dual}), this gives
$\RM (G^\star) = \{ \{x, x\}\mid x\in \N \}$.  Then it is easy to
see that $\mathcal{P}(G^\star) = \{ \{0, x\}\mid x\in \M \}$. Hence, for
the two games in this example we have that $(G^\star)^\star = G$. Neither does
this equality hold in general. (See also Example \ref{Ex2}.)
\end{Ex}

From now onwards we let $k = 2$.

A pair of sequences $(x_n)_{n\in \N}$ and $(y_n)_{n\in \N}$ of
positive integers
is said to be \emph{complementary} if $\{x_n\}\cup \{y_n\} = \N$ and
$\{x_n\}\cap \{y_n\} = \emptyset$.

Let $\alpha < \beta$ be positive irrational numbers
satisfying $1/\alpha + 1/\beta = 1$.
Hence $1 < \alpha < 2 < \beta$. We call $(\alpha, \beta)$ an (ordered)
{\em Beatty pair}. It is well-known \cite{BOHA} that the sequences
$(\lfloor n\alpha \rfloor)_{n\in \N}$ and $(\lfloor n\beta \rfloor)_{n\in \N}$
are complementary.

Our purpose is to prove the following conjecture \cite{DR}:
\begin{Conj}[Duch\^{e}ne-Rigo] Let $(\alpha,\beta)$ be
a Beatty pair. Then there exists an invariant game $G$ such that
$\mathcal{P}(G) = \{\{\lfloor n\alpha \rfloor,
\lfloor n\beta \rfloor \} \mid n \in \M \}$.
\end{Conj}

Let $t\in \N$. We
say that a sequence $(X_n)_{n\in \M}$ of non-negative integers is
{\em $t$-superadditive} if, for all $m, n\in \M$,
\be\label{tsup}
X_m + X_n \leq X_{m + n} <  X_m + X_n + t.
\ee
Note that the left-hand inequality of (\ref{tsup}) is the usual
definition of {\em superadditivity}.

Let $a = (a_n)_{n\in \N}$ and  $b = (b_n)_{n\in \N}$ be sequences
of positive integers and define $a_0 = b_0 = 0$.
We say that the set $\{(a_n, b_n) \mid n\in \M\}$ of ordered pairs is
\emph{$t$-superadditive-complementary},
abbreviated $t$-SAC, if the following criteria are satisfied:
\begin{itemize}
\item $a_1 = 1$,
\item $a$ and $b$ are complementary sequences,
\item $a$ is increasing,
\item $b$ is $t$-superadditive.
\end{itemize}
$\;$ \\
We can now state the main result of this paper :

\begin{Thm}
Suppose that the set
$\{(a_n, b_n) \mid n \in \M\}$ of ordered pairs is $b_1$-SAC.
Define $G$ by setting
$\RM(G) := \{\{a_n, b_n\}\mid n\in \N \}$.
Then
\be\label{result}
\mathcal{P}(G^\star) = \RM(G) \cup \{\boldsymbol 0\}
\ee
and
\be\label{duality}
(G^{\star})^{\star} = G.
\ee
\end{Thm}

\par An immediate consequence of this result is

\begin{Cor}
Suppose that $\{(a_n,b_n) \mid n \in \M \}$ is $b_1$-SAC. Then there is
an invariant game $I$ such that
$\mathcal{P}(I) = \{\{a_n, b_n\} \mid n \in \M \}$.
\end{Cor}

\noindent {\bf Proof of Corollary.} Take $I = G^\star$ in Theorem 1.2.
\hfill $\Box$.\\

It is well-known and easy to check that if $a$ and $b$ are a pair
of complementary homogeneous Beatty sequences, then the set
$\{(a_n, b_n) \mid n \in \M\}$
is $2$-SAC, hence $b_1$-SAC. Therefore, Corollary 1.3 implies Conjecture 1.1.
\par
Because of (\ref{duality}), it is natural to refer to the game
$G^\star$ defined by (1.4)
as the {\em dual} of $G$, when $G$ satisfies the hypotheses of
Theorem 1.2. It is important to note, however,
that the $\lq$duality relation'
(1.7) doesn't always hold for games $G$ not satisfying these hypotheses.

\begin{figure}[ht!]
\centering
\includegraphics[width=0.8\textwidth]{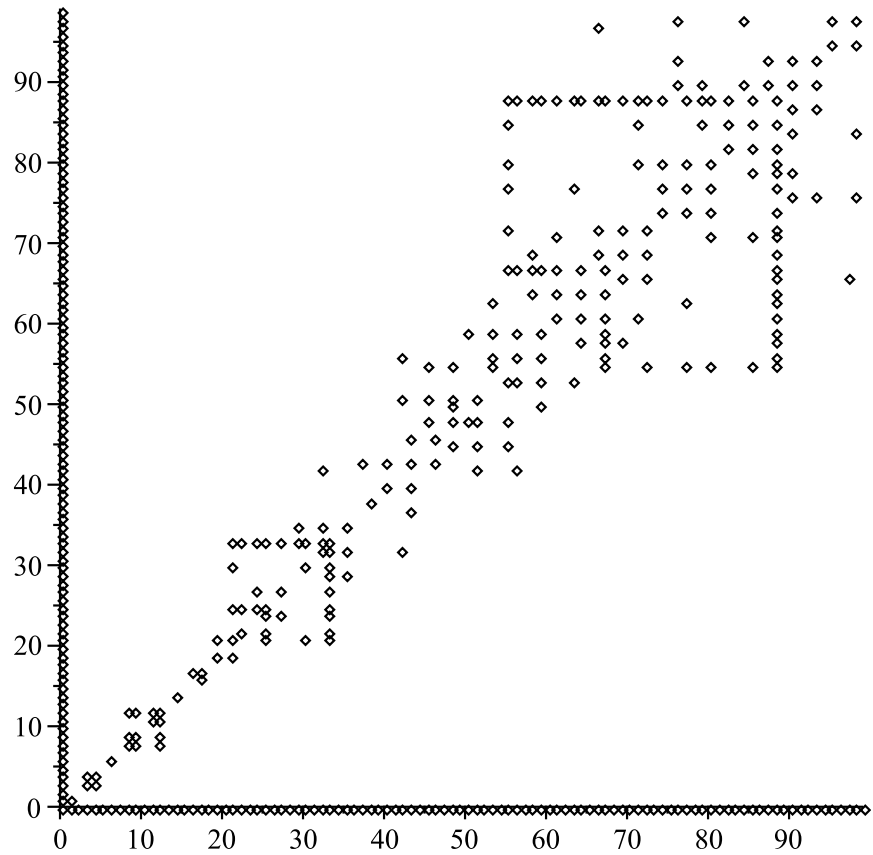}
\includegraphics[width=0.8\textwidth]{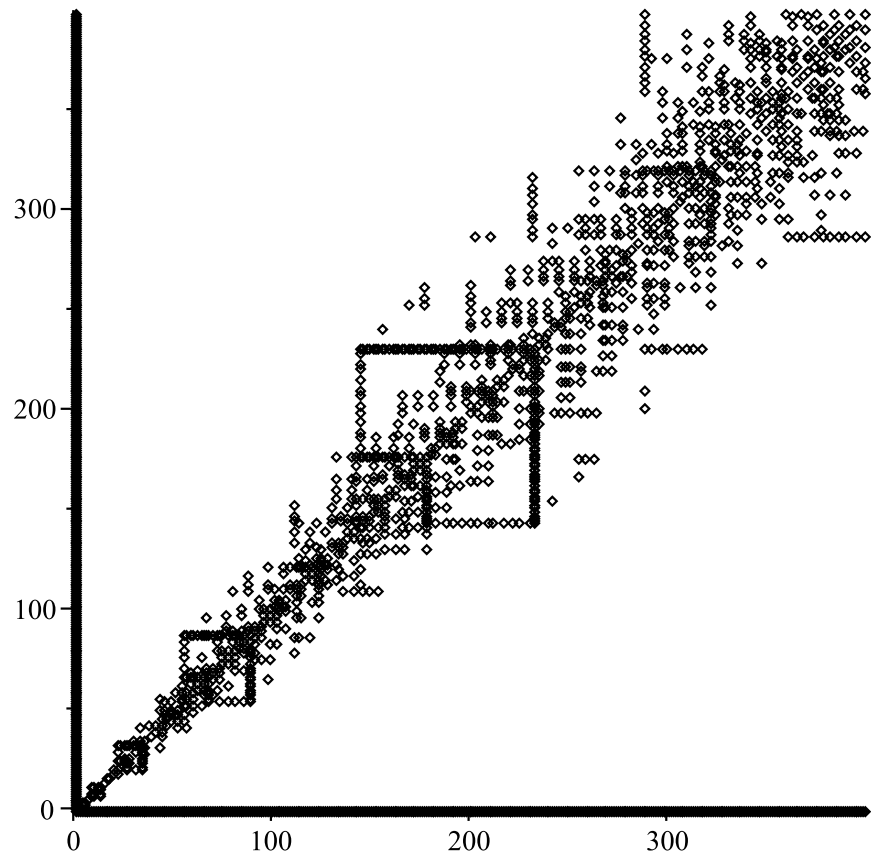}
\caption{The set
$\{ \{i, j\} \in \mathcal{P}(\text{WN}^\star) \mid 0\le i,j\le x \} =
\{\{0, 0\}\}\cup\{ \{i, j\} \in \RM((\text{WN}^\star)^\star)
\mid 0\le i,j\le x \},$ for $x = 100, 400$ respectively.}\label{Fig1}
\end{figure}

\begin{Ex}\label{Ex2}
As a simple but instructive example, take $G =$ WN, the ordinary
Wythoff Nim game \cite{W}, so that
$\RM (\text{WN}) = \{ \{0,i\},(i,i)\mid i\in \N\}$. This set obviously does not
satisfy the hypotheses of Theorem 1.2, whereas $\RM (\text{WN}^\star)$ does so.
Indeed, according to (1.4), we have
\be
\RM ({\hbox{WN}}^{\star}) =
\mathcal{P}({\hbox{WN}})
\backslash \{\boldsymbol 0 \} = \{\{\lfloor n\phi \rfloor,
\lfloor n \phi^2 \rfloor\} \mid n \in \N\}, \;\;\; \phi = \frac{1+\sqrt{5}}{2}.
\ee
It is easy to see that
$\{\{0,x\}\mid x\in \M\} \subset\mathcal{P}({\hbox{WN}}^\star)$. Otherwise
it is easy to check that the first few $P$-positions of ${\hbox{WN}}^\star$ are 
$$(1,1),(3,3),(3,4),(4,4),(6,6),(8,8),(8,9),(8,12),(9,9),(9,12),$$
(see also Figure \ref{Fig1} on page 11) and hence
\begin{eqnarray*}
({\hbox{WN}}^\star)^\star \neq {\hbox{WN}}.
\end{eqnarray*}
But if we go one step further, it follows from (1.8) and Theorem 1.2 that
\begin{eqnarray*}
(({\hbox{WN}}^\star)^\star)^\star = {\hbox{WN}}^\star.
\end{eqnarray*}
In particular, the games WN and (WN$^\star)^\star$ do
have the same $P$-positions.
\end{Ex}

Numerous generalizations and variations of Wythoff Nim
can be found in the literature. In fact, this game can be credited
with opening up the territory of the games we are exploring in this paper.
However, we have not been able to find any literature on the game
$({\hbox{WN}}^\star)^\star$.
\\
\par The rest of the paper is organised as follows. In Section 2,
we will prove Theorem 1.2.
In Section 3,
we explore the problem of describing necessary and sufficient conditions on a
pair $(a_n), (b_n)$
of complementary sequences for there to exist an invariant game
$G$ with $\mathcal{P}(G) = \{\{a_n,b_n\}\} \cup \{\boldsymbol 0 \}$.
We are unable to
solve this problem definitively, though we discuss several pertinent
examples. One of these concerns an application of Theorem 1.2 to
defining an invariant game with
the same solution as the variant game 'the Mouse game' \cite{F3}. In another
example,
we study the set of $P$-positions of the invariant game $G =
(1,2)$-GDWN \cite{L2}.
Here, the $b$-sequence is not increasing and we show that
$\mathcal{P}((G^\star)^\star) \neq \mathcal{P}(G)$.

\setcounter{equation}{0}
\section{Proof of Theorem 1.2}

Let us begin by proving some basic facts about any sequence of $b_1$-SAC
pairs.

\begin{Prop} Suppose that $\{(a_n, b_n) \mid n\in \M\}$ is $b_1$-SAC.
Then, for all $n \in \M$,
\begin{enumerate}[(i)]
\item $b_{n+1} - b_n \geq b_1\ge 2$,
\item $a_{n+1}- a_{n}\in \{1, 2\}$,
\item $a_n < b_n$ and the sequence $(b_n - a_n)$ is non-decreasing,
\item for all $m,n \in \M$,
\be\label{subadd}
a_m + a_n - 1 \leq a_{m+n} \leq a_m + a_n + 1.
\ee
\end{enumerate}
\end{Prop}

\noindent {\bf Proof.} Part (i): By definition $a_1 = 1$. Then,
by complementarity, $b_1\ge 2$. The first inequality
follows by superadditivity.\\
\noindent Part (ii): Let $d_n := a_{n+1} - a_n$. Since $a$ is increasing we
have $d_n \geq 1$ for all $n$. Suppose that there exists an $n$ such that
$d_n\ge 3$. Then, by complementarity,
there exists an $i$ such that $b_i = a_n + 1$ and $b_{i+1} = a_{n} + 2$. But
then $b_{i+1} - b_i = 1$, contradicting (i).\\
\noindent Part (iii): We have $b_1 > a_1$ by definition, and it follows
from parts (i) and (ii) that the sequence $(b_n - a_n)$ is non-decreasing. \\
\noindent Part (iv): Note that, since the sequences $(a_i)$ and $(b_i)$
are increasing and complementary, we have for any $i > 0$ that
\be\label{pineq}
b_{a_i - i} < a_i < b_{a_i - i +1}.
\ee
The inequalities in (\ref{subadd}) are trivial if either $m$ or $n$ equals
zero, so we may suppose that $m, n > 0$. Fix $m$ and $n$.
Let the integers $r, s$ be defined by
\be\label{ineqs}
b_r < a_m < b_{r+1}, \;\;\;\;\;\;
b_s < a_n < b_{s+1}.
\ee
Then, by (\ref{pineq}), it follows that
$a_m = m + r$ and $a_n = n + s$, hence that $a_m + a_n = $ \\ $(m+n) + (r+s)$.
First of all, consider the right-hand inequalities in (\ref{ineqs}).
Superadditivity of $b$ implies that
\begin{eqnarray*}
b_{r+s+2} \geq b_{r+1} + b_{s+1} \geq a_m + a_n + 2 = (m+n) + (r+s+2).
\end{eqnarray*}
Then, by (\ref{pineq}) again we must have
\begin{eqnarray*}
a_{m+n} \leq (m+n) + (r+s+1) = a_m + a_n + 1,
\end{eqnarray*}
which proves the right-hand inequality of (\ref{subadd}).

\par Secondly, the fact that the sequence $b$ is $b_1$-superadditive implies
that
\begin{eqnarray*}
b_{r+s-1} \leq b_{r-1} + b_s + (b_1 - 1) \leq (b_r - b_1) + b_s +
(b_1 - 1) = b_r + b_s - 1.
\end{eqnarray*}
This, together with the left-hand inequalities in (\ref{ineqs}), imply that
\begin{eqnarray*}
b_{r+s-1} \leq (a_m - 1) + (a_n - 1) - 1 = (m+n) + (r+s-3).
\end{eqnarray*}
By complementarity, it follows that
\begin{eqnarray*}
a_{m+n-2} \geq (m+n) + (r+s - 3).
\end{eqnarray*}
Then, the fact that $a$ is increasing implies that
\begin{eqnarray*}
a_{m+n} \geq a_{m+n-2} + 2 \geq (m+n) + (r+s-1) = a_m + a_n - 1,
\end{eqnarray*}
which proves the left-hand inequality of (\ref{subadd}). This
completes the proof of Proposition 2.1.
\hfill $\Box$

\begin{Rem}
In the above proof,
superadditivity of $b$ sufficed, except for the left-hand inequality in (2.1).
Only the latter required $b_1$-superadditivity.
Interestingly enough,
$b_1$-superadditivity is needed for the proof of Theorem 1.2, but
the left-hand inequality in (2.1) is not.
\end{Rem}

For our particular setting,
the next lemma is a special case of part (iii) of the one to follow.
But it is nice to first state it in a more general form.

\begin{Lemma}[A $P$-position is never an invariant move]
Let $G$ be an invariant subtraction game. Then
$\RM (G)\cap \mathcal{P}(G) = \emptyset$.
\end{Lemma}

\noindent{\bf Proof.} Suppose that there was a
move $\boldsymbol r \in \mathcal{P}(G)$. Then, in particular,
$\boldsymbol 0 = \boldsymbol r - \boldsymbol r \in F(\boldsymbol r)$.
But $\boldsymbol 0 \in \mathcal{P}(G)$, so then
$\boldsymbol r \in \mathcal{N}(G)$, a contradiction. \hfill $\Box$\\

The hypothesis of the next lemma is satisfied, in particular, by any
game $G$ for which $\RM(G)\cup\{(\boldsymbol 0)\}$, viewed as an ordered set,
is ($b_1$-)SAC. The items (i) and (ii) characterize precisely the lower ideal
$\mathcal{T}(G)$.

\begin{Lemma}
Let $(a_n)_{n \in \N}$ and $(b_n)_{n \in \N}$ be any pair of
increasing sequences of positive integers, and suppose that
$G$ is an invariant subtraction game with
$\RM (G)= \{\{a_n, b_n \}\}$.
Then
\begin{enumerate}[(i)]
\item $\{0, k\}\in \mathcal{P}(G)$, for all $k \in \M$,
\item $\{k, l\}\in \mathcal{P}(G)$ if $k, l \in \{1, 2, \ldots , b_1 - 1\}$,
\item If $k,l > 0$ then $\{k, l\} \in \mathcal{N}(G)$
if, for some $n > 0$,
\begin{enumerate}[(a)]
\item $k = a_n$ and $l \ge b_n $, or
\item $k = b_n$ and $l \geq a_n$, or
\item $k = a_n$, $a_n = a_{n-1} + 1$ and $b_{n-1} \leq l < b_{n-1} + b_1.$
\end{enumerate}
\end{enumerate}
\end{Lemma}

\noindent{\bf Proof.} Parts (i), (ii):
By the definition of $\RM (G)$, it is clear
that $F(\{k, l \}) = \emptyset $ if
either $\min \{k,l\} = 0$ or $\max \{k,l\} < b_1$.\\
\noindent Part (iii): If (a) holds, then
\begin{eqnarray*}
(k,l) \rightarrow (k,l) \ominus (a_n,b_n) = (0,l-b_n),
\end{eqnarray*}
is an option in $G$. Since $(0,l-b_n) \in \mathcal{P}(G)$ by (i), it
follows that $(k,l) \in \mathcal{N}(G)$. Similarly, if (b) holds then
one considers the option
\begin{eqnarray*}
(k,l) \rightarrow (k,l) \ominus (b_n,a_n) = (0,l-a_n) \in \mathcal{P}(G).
\end{eqnarray*}
Finally, if (c) holds, then we have the option
\begin{eqnarray*}
(k,l) \rightarrow (k,l) \ominus (a_{n-1},b_{n-1}) = (1,l-b_{n-1}).
\end{eqnarray*}
Since $l - b_{n-1} < b_1$, we have $(1,l-b_{n-1}) \in \mathcal{P}(G)$ by
(ii), and hence $(k,l) \in \mathcal{N}(G)$ once more.
\hfill $\Box $

\vspace{2 mm}

\noindent{\bf Proof of Theorem 1.2.} Clearly, (1.7) follows from (1.6) so
it remains to prove the latter. Recall that
the moves in the game $G^\star$ are given by
$\RM(G^\star):=\mathcal{P}(G)\setminus \{\boldsymbol 0\}$ and
where $\RM(G) := \{\{a_n, b_n\} \mid n \in \M\}\setminus \{\boldsymbol 0\}$.
We want to show that
\be\label{RGstar}
\mathcal{P}(G^\star) = \{\{a_n, b_n\} \mid
n \in \M\}.
\ee
By the definition of $\mathcal{P}$, this corresponds to showing that,
for all $(\alpha ,\beta )\in \B$,
\be\label{cond}
\exists \; n \; {\hbox{such that either}} \;
(\alpha,\beta) \rightarrow (a_n,b_n) \; {\hbox{or}} \;
(\alpha,\beta) \rightarrow (b_n,a_n) \; {\hbox{is an option in $G^\star$}}
\ee
if and only if
$\{\alpha ,\beta \} \neq \{a_i, b_i\}$ for all $i \in \M$.
\\
\\
``N$\rightarrow $ P'': Suppose that $\{\alpha,\beta\} \neq \{a_i,b_i\}$ for any
$i \in \M$. If $(\alpha , \beta )\in \mathcal{P}(G)$ then
$(\alpha,\beta) \rightarrow \boldsymbol 0 = (a_0,b_0)$ is
an option in $G^\star$, thus satisfying (\ref{cond}).
If $(\alpha,\beta) \in \mathcal{N}(G)$, then there exists
$(x,y) \in \mathcal{P}(G)$ such that $(\alpha,\beta) \rightarrow
(x,y)$ is an option in $G$. By definition of $\RM (G)$, there exists
$j \in \N$ such that either $(\alpha, \beta) \ominus (a_j,b_j) =
(x,y)$ or $(\alpha,\beta) \ominus (b_j,a_j) = (x,y)$. Note that our
assumptions thus imply that $(x,y) \neq \boldsymbol 0$. Hence $(x,y) \in
\mathcal{P}(G) \backslash \{\boldsymbol 0 \} = \RM (G^\star)$. Since
$(\alpha,\beta) \ominus (x,y) \in \{(a_j,b_j), (b_j,a_j)\}$, we see that
once again (\ref{cond}) is satisfied.
\\
\\
``P$\rightarrow $ N'': Suppose that $\{\alpha,\beta\} = \{a_i,b_i\}$ for
some $i \in \M$ and that (\ref{cond}) holds. By symmetry, it suffices to
consider the following two cases : there exists $m,n \in \M$ such that
$m > 0$ and
either $(a_{m+n},b_{m+n}) \rightarrow (a_n,b_n)$ or $(a_{m+n},b_{m+n})
\rightarrow (b_n,a_n)$ is an option in $G^\star$.
\par First suppose the latter.
Let
\be
(x,y) := (a_{m+n},b_{m+n}) \ominus (b_n,a_n).
\ee
By definition of $G^\star$, we must have $(x,y) \in \mathcal{P}(G)
\backslash \{\boldsymbol 0 \}$. By Lemma 2.2, we may assume that $n > 0$.
Then $x = a_{m+n} - b_n < a_{m+n} - a_n \leq a_m + 1$, by parts (iii) and (iv)
of Proposition 2.1. Hence $x \leq a_m$. By
complementarity, there exists $p \leq m$ such that $x \in \{a_p,b_p\}$.
On the other hand, $y = b_{m+n} - a_n > b_{m+n} - b_n \geq b_m$,
by superadditivity of $b$. In particular, $y > x$. But then $(x,y) \in
\mathcal{N}(G)$, by parts (a),(b) of Lemma 2.3(iii), a contradiction.
\par Second, suppose that $(a_{m+n},b_{m+n}) \rightarrow (a_n,b_n)$ is an
option in $G^\star$. Let
\be\label{diffeq}
(x,y) := (a_{m+n},b_{m+n}) \ominus (a_n,b_n).
\ee
As before, we must prove the contradiction that $(x,y) \in \mathcal{N}(G)$.
By the $b_1$-superadditivity of $b$, we have
\be\label{y}
b_m \leq y < b_m + b_1.
\ee
If $x \leq a_m$ then we can appeal to parts (a),(b) of Lemma 2.3(iii) again.
By the right-hand inequality of (2.1), the only other possibility is that
$x = a_m + 1$. Since $m > 0$ and $y \geq b_m$, we have $y \geq x$. If
$x = b_i$ for some $i$, then part (b) of Lemma 2.3(iii) gives a contradiction.
This leaves the possibility that $x = a_{m+1} = a_m + 1$. But then,
because of (\ref{y}), we get a contradiction from part (c) of Lemma 2.3(iii).
\hfill $\Box$

\setcounter{equation}{0}
\section{Discussion}
In this section we provide four examples and suggest some future work.
\begin{Ex}\label{Ex3}
Let $a$ and $b$ be any complementary, though not necessarily increasing,
sequences beginning as in Table 1 below.
\par As usual, set $a_0 = b_0 := 0$. Note that the set of
pairs $\{(a_n,b_n) \mid n \in
\M\}$ cannot be $b_1$-SAC, since $b_3 = b_{2+1} = b_2 + b_1 + b_1$.
Suppose there were an invariant game $G$ with
$\mathcal{P}(G) = \{\{a_n,b_n\} \mid n \in \M\}$. Then
$(2,6) \in \mathcal{N}(G)$. But $(2,6) = (4,13) \ominus (2,7)$,
a contradiction.

\begin{table}[ht!]
\begin{center}
  \begin{tabular}
{| l || c | c | c |}
    \hline
    $b_n$\T &3&7&13 \\
    $a_n$\T &1&2&4\\ \hline
    $n$  \T &1&2&3\\
    \hline
  \end{tabular}
\end{center}\caption{The $b$-sequence does not satisfy the right-hand
inequality in (1.5).}
\end{table}

\par Nevertheless, if the sequences $a$ and $b$ are increasing, $a_1 = 1$ and
the $b$-sequence grows at only a slightly faster rate than that allowed
by (1.5),
then Theorem 1.2 will hold again. Indeed, suppose that
\be\label{biggerb}
b_2 \geq 2b_1 \;\; {\hbox{and}} \;\; b_{m+n} \geq b_{m+1} + b_n \; {\hbox{for all $m \geq 1$, $n \geq 2$}}.
\ee
We can still use Lemma 2.3 and one may check that the proof of Theorem 1.2 goes
through. Consider (\ref{diffeq}), for example. We still have $x
\leq a_m + 1 \leq a_{m+1}$, since for the right-hand inequality of
(2.1) we only required $b$ to be superadditive. If $n \geq 2$, then
(\ref{biggerb}) implies that $y \geq b_{m+1}$. Then from Lemma 2.3(iii),
parts (a) and (b), it follows that $(x,y) \in \mathcal{N}(G)$.
We can obtain the same conclusion even when $n=1$, since then we still
have $y \geq b_m$ and now $x
= a_{m+1} - 1 < a_{m+1}$, with strict inequality.
\end{Ex}

\begin{Ex}\label{Ex4}
(A similar example to this one appears in \cite{DR}). Let
$a$ and $b$ be any complementary sequences beginning as in Table 2.
\par
Put $a_0 = b_0 := 0$. The set of pairs $\{(a_n,b_n) \mid n \in \M \}$ cannot be
$b_1$-SAC since $b_2 = b_{1+1} = b_1 + b_1 - 1$. Suppose there were an
invariant game $G$ with
$\mathcal{P}(G) = \{\{a_n,b_n\} \mid n \in \M\}$. Then $(1,3) \in
\mathcal{N}(G)$. But $(1,3) = (2,7) \ominus (1,4)$, a contradiction.

\begin{table}[ht!]
\begin{center}
  \begin{tabular}
{| l || c | c |}
    \hline
    $b_n$\T &4&7 \\
    $a_n$\T &1&2\\ \hline
    $n$  \T &1&2\\
    \hline
  \end{tabular}
\end{center}\caption{The $b$-sequence does not satisfy superadditivity,
the left-hand inequality in (1.5).}
\end{table}

This example also arises from a pair of complementary, but
inhomogeneous Beatty sequences.
Let $(\alpha,\beta)$ be a Beatty pair. Let $\gamma, \delta \in
\mathbb{R}$. For each $n \in \N$, let
\be\label{inhom}
a_n := \lfloor n\alpha + \gamma
\rfloor, \;\;\;\;\;\;
b_n := \lfloor n\beta + \delta \rfloor.
\ee
Fraenkel \cite{F1} proved that the sequences $(a_n)$ and $(b_n)$ are
complementary if and only if 
$n\beta + \delta \not\in \mathbb{Z}$ for any 
$n \geq 1$, and 
\be\label{fra}
\frac{\gamma}{\alpha} + \frac{\delta}{\beta} = 0.
\ee
Choose a pair of (small) irrational numbers
$\epsilon_1, \epsilon_2 > 0$. Let $\alpha := \frac{7}{5} + \epsilon_1$,
$\beta := \frac{7}{2} - \epsilon_2$. Choose $\delta \not\in
\mathbb{Q}(\beta)$ satisfying
\be\label{delta}
\frac{1}{2} + \epsilon_2 \leq \delta < 1 - 2\epsilon_2.
\ee
It is not hard to check that, for an appropriate choice of
$\epsilon_1, \epsilon_2, \delta$, the number $\gamma < 0$ defined by
(\ref{fra}) will satisfy
\be\label{gamma}
-\frac{2}{5} - \epsilon_1 \leq \gamma < \frac{1}{5} - 2\epsilon_1.
\ee
From (\ref{delta}) and (\ref{gamma}), one may then verify in turn that
the sequences $(a_n)$ and $(b_n)$ defined by (\ref{inhom}) begin
as in Table 2.
\end{Ex}

\begin{Ex}
For each $n \in \N$, let $a_n := \lfloor \frac{3n}{2} \rfloor$ and
$b_n := 3n-1$. It is easy to see that $(a_n)$ and $(b_n)$ are a pair of
complementary, inhomogeneous Beatty sequences. Put $a_0 = b_0 := 0$, as usual.
In \cite{F3},
a variant game $G$ named $\lq$the Mouse game' was invented with
$\mathcal{P}(G) = \{\{a_n, b_n\} \mid n \in \M\}$.
But, since it is easy to verify that $\{(a_n, b_n)\}$ is $b_1$-SAC,
by Theorem 1.2 we may also introduce an
invariant game $H$, which we call \lq the Mouse trap', with
$\mathcal{P}(H) = \mathcal{P}(G)$. In analogy with Example \ref{Ex2},
the invariant rules are $\mathcal{M}(H) = \mathcal{P}(G^\star).$

\end{Ex}

\begin{Rem}
In \cite{F2, L1}  invariant games with symmetric moves are defined whose
$P$-positions consist of complementary inhomogeneous Beatty sequences (CIBS).
Both papers include variations of Wythoff Nim. In the former
a mis\`ere variation (the player who moves last loses) is studied.
Indeed, we believe it to be the 'most natural/direct'
way to construct a game with CIBS as $P$-positions. In the latter paper,
the terminal positions are $(l, 0)$ and $(0, p-l)$, for some integers
$0<l<p$, so the game is only symmetric if $p = 2l$. Namely, here the game board
is rearranged to
$$\B:= (\M \times \M) \backslash \{(i,j)\mid 0 \le i < l, 0 \le j < p-l \}.$$
\end{Rem}

The above examples provide some extra insight into the
following problem, which nevertheless remains wide open :
\\
\\
{\bf Problem 1.} Let $(a_n)$, $(b_n)$ be a pair of complementary,
increasing sequences with $a_1 = 1$. Find necessary and sufficient conditions
for the existence of an invariant game $G$ with $\mathcal{P}(G) =
\{\{a_n,b_n\}\} \cup \{\boldsymbol 0 \}$.
\\
\par
A special case which might be more tractable is the case of
inhomogeneous Beatty sequences. Motivated by Examples 4
and 5, we may ask
\\
\\
{\bf Problem 2.} Let $(a_n)$, $(b_n)$ be a pair of complementary,
inhomogeneous Beatty sequences with $a_1 = 1$. Is it true that there
exists an invariant game $G$ with $\mathcal{P}(G) =
\{\{a_n,b_n\}\} \cup \{\boldsymbol 0 \}$ if and only if the set of
pairs $\{(a_n,b_n)\}$ is $b_1$-SAC ?
\\
\par
In studying these problems, it is natural to ask whether the method of
Theorem 1.2 will ever fail, in the following sense :
\\
\\
{\bf Problem 3.} Does there exist a pair $(a_n)$, $(b_n)$ of
complementary, increasing sequences, with $a_1 = 1$, such that
$\{\{a_n,b_n\}\} \cup \{\boldsymbol 0 \} = \mathcal{P}(G)$ for some
invariant game $G$, but $\mathcal{P}((G^\star)^\star) \neq \mathcal{P}(G)$ ?
\\
\par
We know that the answer to Problem 3 is yes,
if we drop the condition that the sequences
$(a_n)$, $(b_n)$ be increasing. Consider the following example :

\begin{Ex}
Let $G$ be the invariant game $(1,2)$-GDWN, studied in \cite{L2},
so that $\RM (G) = \{ \{0,i\},(i,i),\{i,2i\}\mid i\in \N\}.$
Define
\begin{eqnarray*}
\{\{a_n, b_n\} \mid n \in \M\} := \mathcal{P}(G), \;\;\;
{\hbox{where $(a_n)$ is increasing}}.
\end{eqnarray*}
Then the sequences $(a_n)_{n \in \N}$ and
$(b_n)_{n \in \N}$ are complementary, but $b$ is not increasing. Table 3 gives
the initial $P$-positions of this game.

\begin{table}[ht!]
\begin{center}
  \begin{tabular}
{| l || c | c | c | c | c | c | c | c | c | c | c |c|}
    \hline
    $b_n$\T &0&3&6&5&10&14&17&25&28&18&35&23 \\
    $a_n$\T &0&1&2&4&7&8&9&11&12&13&15&16\\ \hline
    $n$  \T &0&1&2&3&4&5&6&7&8&9&10&11\\
    \hline
  \end{tabular}
\end{center}\caption{The initial $P$-positions of the game $(1,2)$-GDWN.
For a more comprehensive list, see \cite{L2}.}
\end{table}

Now consider the game $G^\star$, as defined by (1.4). It is not hard to check
that $(11,23) \in \mathcal{P}(G^\star)$.
However, by brute-force calculation one
may also verify that $(104,235)$ and $(115,258)$ are in $\mathcal{P}(G)$.
Since
\be
(115,258) \ominus (104,235) = (11,23),
\ee
we see that $\mathcal{P}((G^\star)^\star)$ cannot coincide with $\mathcal{P}(G)$.
\end{Ex}

Another possible direction for future work is to extend our
results in some manner to $k$-pile subtraction games for $k > 2$, or
even perhaps to consider subtraction games played on other
partially-ordered semigroups. Alternatively, one might try to
extend the notion of $\lq$invariance' to games which cannot be
formulated as subtraction games. Many such games appear in the
literature, see for example \cite{S}, where 14 such games are proved 
Pspace-complete, 3 played on graphs, including Geography --- whose many 
variations have  been addressed in other papers --- and 11 on 
propositional formulas. Another example is annihilation games --- if a 
token moves onto another one, both disappear --- for which there is 
a polynomial-time winning strategy \cite{FY}.

Finally, the ``$\star$-operator'' introduced
in (1.4) and the duality in (1.7)
may turn out to be useful in other contexts.

\section*{Acknowledgements}

Joint work on this project was initiated while
the first and second authors were guests at the
Weizmann Institute of Science, and we thank our hosts for their
hospitality. The research of the
second author is supported by a grant from
the Swedish Science Research Council (Vetenskapsr\aa det).

\end{document}